# Functions on space curves

V. Goryunov [*]

**Abstract**

We classify simple singularities of functions on space curves. We show that their bifurcation sets have properties very similar to those of functions on smooth manifolds and complete intersections [1, 2]: the $k(\pi, 1)$-theorem for the bifurcation diagram of functions is true, and both this diagram and the discriminant are Saito's free divisors.

## 1 Space curves

A germ of any reduced space curve is a determinantal variety: it is a zero set of all order $n$ minors of a germ of some $n \times (n+1)$-matrix $M$ on $\mathbf{C}^3$. In what follows we often treat a space curve as such a matrix.

**Definition 1.1** *Two curve-germs at the origin, given by matrices $M$ and $M'$, are said to be* equivalent *if there exist two matrix-germs $A$ and $B$ and a diffeomorphism-germ $h$ of $(\mathbf{C}^3, 0)$ such that*

$$AMB = M' \circ h \ .$$

We say that *corank* of a space curve-germ is $c$ if the rank at the origin of its defining matrix is $n - c$. Obviously, such a curve can be given by a $c \times (c+1)$-matrix which is just a zero matrix at the origin.

Classification of simple space curves within this approach has been obtained in [8]. As a subset it contains simple plane $A, D, E$-curves and Giusti's list of simple 1-dimensional complete intersections in the 3-space [9, 10, 2]. Those are corank 1 curves. All other simple curves have corank 2. Since we will not need the entire list we are not reproducing it here.

---

[*]Partially supported by the grants from The Royal Society and The Fields Institute for Research in Mathematical Sciences.

# 2  Simple function singularities

**Definition 2.1** *A function on a space curve is a pair of germs $(M, f)$ on $\mathbf{C}^3$, where $M$ is an $n \times (n+1)$-matrix and $f$ a function.*

**Definition 2.2** *Two germs of functions on curves at the origin, $(M, f)$ and $(M', f')$, are $\mathcal{R}_c$-equivalent if there exist two matrix-germs $A$ and $B$, a diffeomorphism-germ $h$ of $(\mathbf{C}^3, 0)$ and a function $g$ in the ideal generated by the maximal minors of $M$ such that*

$$(AMB, f + g) = (M' \circ h, f' \circ h) \ .$$

Notation $\mathcal{R}_c$ is chosen to indicate that this is a sort of right equivalence of functions, but this time restricted to curves.

The notion of $\mathcal{R}_c$-equivalence satisfies all the conditions of Damon's good geometrical equivalence [6, 7]. Thus all standard theorems like those of versality and finite determinacy are valid in our case. We can apply traditional technique [3, 1] to classify singularities and write out their $\mathcal{R}_c$-versal deformations.

One of our aims is to classify $\mathcal{R}_c$-simple singularities. This implies that the participating curve has to be simple as a determinantal variety. Thus we have to study functions either on simple plane curves, or on non-planar 1-dimensional complete intersections in $\mathbf{C}^3$, or on determinantal space curves of corank 2. The second case is easily seen to contain no $\mathcal{R}_c$-simple functions. The simple lists for the two others are given below.

## 2.1  Functions on plane curves

**Theorem 2.3** [11] *The complete list of $\mathcal{R}_c$-simple functions on plane curves as follows:*

| notation | curve equation | function | restrictions |
|---|---|---|---|
| $A_k$ | $y$ | $x^{k+1}$ | $k \geq 0$ |
| $C_{p,q}$ | $xy$ | $x^p + y^q$ | $p \geq q \geq 1$ |
| $B_k$ | $x^2 + y^k$ | $y$ | $k \geq 3$ |
| $F_{2k+1}$ | $x^2 + y^3$ | $y^k$ | $k \geq 2$ |
| $F_{2k+4}$ | $x^2 + y^3$ | $xy^k$ | $k \geq 0$ |

All the adjacencies of the singularities are compositions of those within the series obtained by reducing the indices (we additionally set $B_2 = C_{1,1}$, $F_3 = B_3$) and
$$F_{p+q+1} \to C_{p,q} \to A_{p+q-1} \ .$$

There are just two bounding functions:
$$\begin{array}{lll} X_9^* : & x + \alpha y^2 & \text{on} \quad x^2 + y^4 = 0 \\ J_{10}^* : & x + \alpha y & \text{on} \quad x^3 + y^3 = 0 \ . \end{array}$$

Here $\alpha$ is a generic complex number (modulus). Both singularities are adjacent to $B_4$ and $F_4$. Also $X_9^*$ is adjacent to $C_{3,1}$.

## 2.2 Functions on space curves

**Theorem 2.4** *The list of simple functions on determinantal space curves of corank 2 is as follows:*

| notation | curve matrix | function | restrictions |
|---|---|---|---|
| $C_{p,q,r}$ | $\begin{vmatrix} x & y & 0 \\ 0 & y & z \end{vmatrix}$ | $x^p + y^q + z^r$ | $p \geq q \geq r \geq 1$ |
| $\dot{F}_{2k+2}$ | $\begin{vmatrix} x & y & 0 \\ y^2 & x & z \end{vmatrix}$ | $z + y^k$ | $k \geq 1$ |
| $\dot{F}_{2k+5}$ | $\begin{vmatrix} x & y & 0 \\ y^2 & x & z \end{vmatrix}$ | $z + xy^k$ | $k \geq 0$ |
| $\check{E}_6$ | $\begin{vmatrix} x & y & z \\ z^2 & x & y \end{vmatrix}$ | $z$ | — |
| $\check{E}_7$ | $\begin{vmatrix} x & y & z \\ yz & x & y \end{vmatrix}$ | $z$ | — |
| $\check{E}_8$ | $\begin{vmatrix} x & y & z \\ z^3 & x & y \end{vmatrix}$ | $z$ | — |

Since the classification technique is very similar to that used in numerous earlier classifications there is no reason to give the proof of this theorem here.

One obtains obvious adjacencies of the singularities reducing the indices in the series. Some other adjacencies showing the relations between the series are:

$$C_{p,q,r} \to C_{p,q} \qquad \dot{F}_r \to \dot{F}_{r-1} \qquad \dot{F}_{p+q+3} \to C_{p,q,1} \qquad \check{E}_6 \to \dot{F}_5 \ .$$

The two bounding singularities from the previous subsection serve as a complete bounding list in this case as well.

## 3 Properties of functions on curves

### 3.1 Tjurina number and Milnor number

Let $(x_1, x_2, x_3)$ be coordinates on $(\mathbf{C}^3, 0)$ and $\mathcal{O}_3$ the space of holomorphic function-germs on it. Choosing an ordering of entries of pairs consisting of an $n \times (n+1)$-matrix and a function, we identify the space of their germs with the module $\mathcal{O}_3^{n(n+1)+1}$. The tangent space $T(M, f)$ to the (extended) $\mathcal{R}_c$-equivalence class of a germ $(M, f)$ in this module is the $\mathcal{O}_3$-module generated by the elements

$$\begin{aligned}
(E_{ij}^n M, 0), & \qquad i,j = 1, \ldots, n, \\
(M E_{kl}^{n+1}, 0), & \qquad k,l = 1, \ldots, n+1, \\
(0, \partial f / \partial x_m), & \qquad m = 1, 2, 3.
\end{aligned}$$

Here $E_{ij}^n$ is the $n \times n$-matrix having 1 in the $ij$th place and zeros everywhere else.

We set $\tau(M, f)$ to be the dimension of the linear space $\mathcal{O}_3^{n(n+1)+1}/T(M, f)$ and call it the *Tjurina number* of the singularity. This is the dimension of the base of an $\mathcal{R}_c$-miniversal deformation of $(M, f)$. Such a deformation can be taken in the form

$$(M, f) + \sum_{i=0}^{\tau-1} \lambda_i e_i \ ,$$

where the $e_i$ are elements of $\mathcal{O}_3^{n(n+1)+1}$ which project to a linear basis of the quotient $\mathcal{O}_3^{n(n+1)+1}/T(M,f)$.

For singularities $C_{p,q}$ and $C_{p,q,r}$ we have $\tau = p+q$ and $\tau = p+q+r+1$ respectively. For all the other simple singularities of our lists, $\tau$ is the subscript in the notation. Also $\tau(X_9^*) = \tau(J_{10}^*) = 6$.

Another characteristic of a function-germ on a curve is the number $\mu(M,f)$ (*Milnor number*) of Morse critical points which a generic small perturbation of $f$ has on a generic smoothing of curve $M$.

**Conjecture 3.1** $\quad \tau(M,f) = \mu(M,f)$ .

The conjecture is true for functions on complete intersections [11], for $\mathcal{R}_c$-simple singularities and in all the examples the author knows.

## 3.2 Discriminant as a free divisor

**Definition 3.2** *Consider the base $\mathbf{C}^\tau$ of an $\mathcal{R}_c$-miniversal deformation of singularity $(M,f)$. The discriminant $\Delta(M,f) \subset \mathbf{C}^\tau$ of $(M,f)$ is the closure of the set of those values of the deformation parameters for which the corresponding curve is smooth and the function on it has a critical point on its zero level.*

We recall that a hypersurface $\Gamma$ in $N$-dimensional complex linear space is called *a free divisor* if the algebra $\Theta_\Gamma$ of vector fields on $\mathbf{C}^N$ tangent to it (that is preserving its ideal) is generated by $N$ elements as a module over functions on $\mathbf{C}^N$.

**Theorem 3.3** *Assume $\tau(M,f) = \mu(M,f)$. Then the discriminant $\Delta(M,f) \subset \mathbf{C}^\tau$ is a free divisor.*

*Proof* (cf. [16, 17, 18, 12, 13, 2]). Let $(\mathcal{M}, F) = (\mathcal{M}(x,\lambda), F(x,\lambda))$ be an $\mathcal{R}_c$-miniversal deformation of $(M,f)$, with $\lambda = (\lambda_0, \ldots, \lambda_{\tau-1}) \in \mathbf{C}^\tau$ being the parameters. For any $i = 0, \ldots, \tau-1$, due to the versality there exists a decomposition

$$F\frac{\partial}{\partial \lambda_i}(\mathcal{M}, F) = (\mathcal{AMB}, \mathcal{G}) + \sum_{r=1}^{3} h_{ir}\frac{\partial}{\partial x_r}(\mathcal{M}, F) + \sum_{j=0}^{\tau-1} v_{ij}\frac{\partial}{\partial \lambda_j}(\mathcal{M}, F) ,$$

where $\mathcal{A}(x,\lambda)$ and $\mathcal{B}(x,\lambda)$ are matrix-germs, $h_{ir}(x,\lambda)$ and $v_{ij}(\lambda)$ function-germs, and $\mathcal{G}(x,\lambda)$ is an element of the ideal generated by the maximal minors of $\mathcal{M}$. A functional factor or differentiation in front of a pair (matrix, function) means multiplication by the function or corresponding differentiation of both items.

The vector fields $\nu_i = \sum_{j=0}^{\tau-1} v_{ij}(\lambda)\partial_{\lambda_j}$ are tangent to $\Delta$. Assume the deformation $(\mathcal{M}, F)$ is monomial and $\lambda_0$ is the free term of $F$. Then it is easily verified that $\det(v_{ij}(\lambda_0, 0, \ldots, 0)) = \lambda_0^\tau$. On the other hand, $\det(v_{ij})$ has to vanish on $\Delta$, and hence is proportional to its defining equation. Since the latter is a polynomial of degree $\mu = \tau$ in $\lambda_0$, $\det(v_{ij}) = 0$ may be taken to be such an equation itself. This implies that the vector fields $\nu_0, \ldots, \nu_{\tau-1}$ generate $\Theta_\Gamma$ as a free module over $\mathcal{O}_\tau$. $\square$

Thus, for example, the discriminant of any $\mathcal{R}_c$-simple singularity is a free divisor. Modulo Conjecture 3.1 this is true for any singularity of finite $\mathcal{R}_c$-codimension.

## 3.3 Bifurcation diagram of functions as a free divisor

Consider a *trunkated* $\mathcal{R}_c$-miniversal deformation $(\mathcal{M}', F')$ of a singularity $(M, f)$, that is one allowing just functions vanishing at $0 \in \mathbf{C}^3$. Its base is of dimension $\tau(M, f) - 1$. Note that in this case the deformation $(\mathcal{M}', F' + \lambda_0)$, where $\lambda_0$ is an additional parameter, is $\mathcal{R}_c$-miniversal for $(M, f)$.

**Definition 3.4** *The bifurcation diagram of functions $\Sigma(M, f) \subset \mathbf{C}^{\tau-1}$ is the set of those values of parameters of the trunkated deformation for which either the corresponding curve is not smooth or the function on it has either a degenerate critical point or at least two critical points on the same level.*

In general, $\Sigma$ has three irreducible components responsible for the degenerations mentioned in the definition.

**Theorem 3.5** *Assume $\tau(M, f) = \mu(M, f)$. Then the bifurcation diagram of functions $\Sigma(M, f) \subset \mathbf{C}^{\tau-1}$ is a free divisor.*

A proof of this statement is absolutely similar to those for functions on smooth manifolds [5, 15] and for functions on complete intersections [12, 13].

The generators of $\Theta_\Sigma$

$$\omega_i = \sum_{j=1}^{\tau-1} w_{ij}(\lambda')\partial_{\lambda_j}, \qquad i = 1,\ldots,\tau-1, \qquad \lambda' = (\lambda_1,\ldots,\lambda_{\tau-1}) \in \mathbf{C}^{\tau-1},$$

are obtained from the decompositions

$$F^i \tfrac{\partial}{\partial \lambda_i}(\mathcal{M}', F') =$$

$$= (\mathcal{A}'\mathcal{M}'\mathcal{B}', \mathcal{G}') + \sum_{r=1}^{3} h'_{ir} \frac{\partial}{\partial x_r}(\mathcal{M}', F') + w_{i0} + \sum_{j=1}^{\tau-1} w_{ij}\frac{\partial}{\partial \lambda_j}(\mathcal{M}', F'),$$

where $\mathcal{A}'(x,\lambda')$ and $\mathcal{B}'(x,\lambda')$ are matrix-germs, $h'_{ir}(x,\lambda')$ and $w_{ij}(\lambda')$ function-germs, and $\mathcal{G}'(x,\lambda')$ is an element of the ideal generated by the maximal minors of $\mathcal{M}'$.

**Example 3.6** For the singularity $C_{1,1,1}$ one can take a trunkated $\mathcal{R}_c$-mini-versal deformation in the form

$$\left( \begin{vmatrix} x & y & \alpha \\ \beta & y+\gamma & z \end{vmatrix}, \quad x+y+z \right).$$

The algebra of vector fields on $\mathbf{C}^3$ tangent to the bifurcation diagram $\Sigma(C_{1,1,1})$ is a free module over $\mathcal{O}_3$ generated by the fields of degree 1 (the Euler field), 2 and 3.

The projectivisation of $\Sigma(C_{1,1,1})$ in $\mathbf{C}P^2$ is a nodal cubic with its three tangents lines at the inflection points. The cubic corresponds to degenerate critical points on smooth curves and the three lines to non-smooth curves. On the left-hand side of the picture below we show this projectivisation (the isolated point in the centre is the node of the cubic). On the right-hand side of the same picture there is given the bifurcation diagram in $\mathbf{R}^3$ for the other obvious version of the complex singularity $C_{1,1,1}$: inside the cone bounded by the cubic there is situated the isolated real line of intersection of the two planes tangent to the cubic surface along its complex parabolic lines.

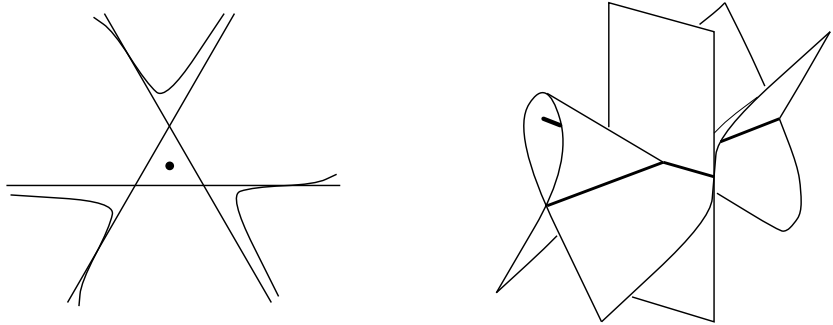

## 3.4 Lyashko-Looijenga mapping

Let $\mathbf{C}^{\tau-1}$ be the space of all monic polynomials in one variable of degree $\tau$ with vanishing sum of the roots, and $\Xi \subset \mathbf{C}^{\tau-1}$ the set of polynomials with multiple roots.

For a singularity $(M, f)$, with Tjurina number $\tau$ equal to its Milnor number, consider the mapping

$$\mathbf{C}^{\tau-1} \setminus \Sigma(M, f) \to \mathbf{C}^{\tau-1} \setminus \Xi$$

from the complement of the bifurcation diagram of functions which sends a Morse function on a smooth curve to the unordered set of its critical values shifted by their arithmetic mean. This mapping is easily verified to be extendible to that between the ambient complex linear spaces [11]. We call the extension *the Lyashko-Looijenga mapping* and denote it by $LL$.

**Theorem 3.7** *For an $\mathcal{R}_c$-simple function on a space curve, the Lyashko-Looijenga mapping is a finite covering. As a mapping from $\mathbf{C}^{\tau-1} \setminus \Sigma$ to $\mathbf{C}^{\tau-1} \setminus \Xi$ it has no branching.*

Since the space $\mathbf{C}^{\tau-1} \setminus \Xi$ is a classification space of the Artin's braid group $B(\tau)$ on $\tau$ threads, we get

**Corollary 3.8** *For an $\mathcal{R}_c$-simple singularity $(M, f)$, the complement to its bifurcation diagram of functions in the base of its truncated $\mathcal{R}_c$-miniversal deformation is a $k(\pi, 1)$-space, where $\pi$ is a subgroup of finite index in the group of braids on $\tau(M, f)$ threads.*

The index $[B(\tau) : \pi]$ is the degree of the mapping $LL : \mathbf{C}^{\tau-1} \to \mathbf{C}^{\tau-1}$. The latter is easy to find taking the truncated deformation of our table normal form to be monomial and hence quasi-homogeneous, and recalling that the degree of a finite quasi-homogeneous mapping is the ratio of the products of the weights of the coordinate functions and of the arguments. Thus obtained indices are given below:

| $A_k$ | $C_{p,q}$ | $F_r$ | $B_k$ |
|---|---|---|---|
| $(k+1)^{k-1}$ | $\frac{(p+q-1)!\,p^p\,q^q}{(p-1)!\,(q-1)!}$ | $\frac{(r-2)(r-1)^r r}{24}$ | $1$ |

| $C_{p,q,r}$ | $\dot{F}_r$ | $\check{E}_6$ | $\check{E}_7$ | $\check{E}_8$ |
|---|---|---|---|---|
| $\frac{(p+q+r+1)!\,p^p\,q^q\,r^r}{(p-1)!\,(q-1)!\,(r-1)!}$ | $\frac{(r-3)^r(r-2)(r-1)r}{24}$ | $3^5$ | $2^7 \cdot 7$ | $2^4 \cdot 3^5$ |

The Theorem and its Corollary are analogous to the classical theorem on simple functions on smooth manifolds [14, 4, 1] and generalise a similar assertion about simple functions on planar curves [11, 2]. The proof is absolutely similar to those cases and we work it out in details just for one series which does not require too long calculations.

*Proof of Theorem 3.7 for singularities $C_{p,q,r}$.* It is more convenient to prove the theorem for a cylindrical version of the Lyashko-Looijenga mapping defined on the base $\mathbf{C}^\tau$ of (non-truncated) $\mathcal{R}_c$-miniversal deformation of a singularity. It sends a generic point of $\mathbf{C}^\tau$ to the unordered set of critical values of the corresponding function without any shifting of them, that is to a point of another copy of $\mathbf{C}^\tau$. The mapping uniquely extends to non-generic points as well. It will be denoted by $\widetilde{LL}$. In fact this is just a trivial one-parameter unfolding of the mapping introduced initially.

Now consider an $\mathcal{R}_c$-miniversal deformation of singularity $C_{p,q,r}$:

$$\begin{vmatrix} x & y & \alpha \\ \beta & y+\gamma & z \end{vmatrix},$$

$$x^p + \lambda_{1,1}x^{p-1} + \ldots + \lambda_{1,p-1}x + y^q + \lambda_{2,1}y^{q-1} + \ldots + \lambda_{2,q-1}y +$$

$$+ z^r + \lambda_{3,1}x^{r-1} + \ldots + \lambda_{3,r-1}z + \lambda_0 ,$$

with all the variables, except for $x, y, z$, being the parameters.

To prove the theorem is to show that

- the Lyashko-Looijenga mapping is a local diffeomorphism out of the bifurcation diagram of functions;
- $\widetilde{LL}^{-1}(0) = 0$.

Let us prove the first of the claims. Assuming the above determinantal curve is not singular we express $x$ and $z$ in terms of $y$ and obtain a function in just one variable $y$:

$$F(y) = (\tfrac{\beta y}{y+\gamma})^p + \lambda_{1,1}(\tfrac{\beta y}{y+\gamma})^{p-1} + \ldots + \lambda_{1,p-1}\tfrac{\beta y}{y+\gamma} +$$

$$+ y^q + \lambda_{2,1} y^{q-1} + \ldots + \lambda_{2,q-1} y +$$

$$+ (\tfrac{\alpha(y+\gamma)}{y})^r + \lambda_{3,1}(\tfrac{\alpha(y+\gamma)}{y})^{r-1} + \ldots + \lambda_{3,r-1}\tfrac{\alpha(y+\gamma)}{y} + \lambda_0 \ .$$

We have to demonstrate that evaluating the derivatives of this function with respect to all its $p+q+r+1$ parameters at all the $p+q+r+1$ critical points of $F$, provides a non-degenerate matrix. It is sufficient to show that the system of the functions

$$\tfrac{\partial F}{\partial \beta}, (\tfrac{\beta y}{y+\gamma})^{p-1}, \ldots, \tfrac{\beta y}{y+\gamma}, y^{q-1}, \ldots, y, \tfrac{\partial F}{\partial \alpha}, (\tfrac{\alpha(y+\gamma)}{y})^{r-1}, \ldots, \tfrac{\alpha(y+\gamma)}{y}, 1, \tfrac{\partial F}{\partial \gamma}$$

is linearly independent at $p+q+r+1$ distinct arbitrary points $y_i \neq 0, -\gamma$ on the $y$-axis.

Since $\alpha\beta\gamma \neq 0$, this is to similarly test the system

$$(\tfrac{y}{y+\gamma})^p, (\tfrac{y}{y+\gamma})^{p-1}, \ldots, \tfrac{y}{y+\gamma}, y^{q-1}, \ldots, y,$$

$$(\tfrac{y+\gamma}{y})^r, (\tfrac{y+\gamma}{y})^{r-1}, \ldots, \tfrac{y+\gamma}{y}, 1, \tfrac{p(\beta y)^p}{(y+\gamma)^{p+1}} + \tfrac{r(\alpha(y+\gamma))^r}{y^{r+1}}$$

or, equivalently, the system

$$(y+\gamma)^{-p}, \ldots, (y+\gamma)^{-1}, y^{q-1}, \ldots, y, 1, y^{-1}, \ldots, y^{-r}, A(y+\gamma)^{-p-1} + By^{-r-1},$$

where $A$ and $B$ are certain non-zero constants.

Multiplication of the last system by $y^{r+1}(y+\gamma)^{p+1}$ transforms it into

$$y^{r+1}(y+\gamma), \ldots, y^{r+1}(y+\gamma)^p, y^{q+r}(y+\gamma)^{p+1}, \ldots, y(y+\gamma)^{p+1},$$

$$Ay^{r+1} + B(y+\gamma)^{p+1}.$$

Add an extra function $y^{r+1}$ to this system and an extra point $y_{p+q+r+2} = 0$ to the original $(p+q+r+1)$-tuple of points $y_i$ of evaluation. Thus obtained square evaluation matrix of order $p+q+r+2$ is easily seen to be non-degenerate. On the other hand, its determinant is the determinant of our original evaluation matrix of order $p+q+r+1$ multiplied by $B\gamma^{p+1}$.

Now we check that $\widetilde{LL}^{-1} = \{0\}$, that is the only member of the above miniversal family with the only critical value 0 (hence, of multiplicity $\tau = p+q+r+1$) is the undeformed singularity. Thus let us search for such a member in all possible different circumstances.

First assume that in our miniversal family we have a function on a smooth curve with only one critical value $a$. Let $k$ be the number of its geometrically distinct critical points. Those are the roots of the numerator of the rational equation $F(y) - a = 0$, where $F$ is the function considered above. The total multiplicity of these roots has to be $p+q+r+1+k$. But the degree of the numerator is just $p+q+r$.

Now let the curve have one node. Then this is a hyperbola $C$ with a straight line $\ell$ intersecting it at some point $m$ and parallel to one of the coordinate axes. Assuming $\ell$ to be of the $x$-direction we get $C$ situated in the coordinate $yz$-plane. According to the construction of the Lyashko-Looijenga mapping as an extension of the mapping from the complement to the bifurcation diagram of functions, the value of a generic function at a node is a critical value (of multiplicity 2). Thus, if we want our function to have just one critical value (in particular, on $\ell$), it has to contain only monomial $x^p$ in $x$. This provides contribution $p$ from $\ell$ to the total multiplicity of the critical value.

The restriction $\phi$ of our function to $C$ has all its critical points on the level of point $m$. Write the equation $\phi(y,z) - \phi(m) = 0$ on $C$ as a rational equation in one variable. The degree of the numerator of this equation is $q+r$. If, besides $m$, $\phi$ has $k$ distinct critical points on $C$, the total contribution from $C$ to the multiplicity of our only critical value is $q+r-k$. Recalling the contribution from the line $\ell$ we see that the multiplicity of this critical value $\phi(m)$ on the entire curve is at most $p + q + r$, that is 1 less than would be required for the corresponding point of the base space to be in $\widetilde{LL}^{-1}(0)$.

Similar considerations in the case of two nodes (when the curve is just 3 straight lines) show that again the only critical value could have multiplicity at most $p + q + r$.

Finally, let the curve be the most degenerate, that is consisting of the three coordinate axes. Restriction of the function to each of them has to have just one critical point, namely, the origin. Since we want this value to be zero, the only way to achieve this is to take the undeformed singularity.

**Acknowledgements.** The author is very thankful to The Fields Institute for Research in Mathematical Sciences for its kind hospitality and nice working atmosphere during the visit in June 1997 when this work was completed.

Division of Pure Mathematics
Department of Mathematical Sciences
The University of Liverpool
Liverpool L69 3BX
UK

e-mail: `goryunov@liv.ac.uk`